\documentclass[11pt]{article}
\usepackage{amsthm,epsfig,a4}
\usepackage[ansinew]{inputenc}
\usepackage{amsmath,amssymb}
\usepackage{amsfonts}
\usepackage{eufrak}
\usepackage{diagrams}
\def\q{\nolinebreak \hfill $\Box$} %q.e.d.%
\def\N{\\[5mm] \indent} %Beweisende Neuer Absatz%

\def\:{\colon\thinspace}
\theoremstyle{plain}
\newtheorem{Lemma}{Lemma}[section]
\newtheorem{Prop}[Lemma]{Proposition}
\newtheorem{Thm}[Lemma]{Theorem}

\theoremstyle{definition}
\newtheorem{Def}[Lemma]{Definition}

\newtheorem*{Ex}{Example}
\newtheorem*{Rem}{Remark}
\newtheorem*{Ac}{Acknowledgement}
\allowdisplaybreaks
\begin{document}
\title{On rational homotopy and minimal models}
\author{\textsc{Christoph Bock}
%\\ \small{e-mail: bock@mi.uni-erlangen.de}
}
\date{}
\maketitle

\begin{abstract}
We prove a result that enables us to calculate the rational homotopy of a wide class of spaces by the theory of minimal models. 
%The latter are derivided from properties of the de Rham complex.
\end{abstract}
\section{Introduction} \label{Intro}
Let $\mathbb{K}$ be a field of characteristic zero. A
\emph{differential graded algebra (DGA)}\index{Differential Graded
Algebra (DGA)} is a graded $\mathbb{K}$-algebra $A = \bigoplus_{i
\in \mathbb{N}}~ A^i$ together with a $\mathbb{K}$-linear map $d \:
A \to A$ such that $d(A^i) \subset A^{i+1}$ and the following
conditions are satisfied:
\begin{itemize}
\item[(i)] The $\mathbb{K}$-algebra structure of $A$ is given by an
inclusion $\mathbb{K} \hookrightarrow A^0$.

\item[(ii)] The multiplication is graded commutative, i.e.\ for $a \in
A^i$ and $b \in A^j$ one has $a \cdot b = (-1)^{i \cdot j} b \cdot a
\in A^{i+j}$.

\item[(iii)] The Leibniz rule holds: $\forall_{a \in A^i} \forall_{b
\in A} ~ d(a \cdot b) = d(a) \cdot b + (-1)^i a \cdot d(b)$

\item[(iv)] The map $d$ is a differential, i.e.\ $d^2 = 0$.
\end{itemize}
Further, we define $|a| := i$ for $a \in A^i$.

\begin{Ex}
Given a manifold $M$, one can consider the
complex of its differential forms $(\Omega(M),d)$, which has the
structure of a differential graded algebra over the field $\mathbb{R}$.
\end{Ex}

The \emph{$i$-th cohomology of a DGA} $(A,d)$ is the algebra
$$ H^i(A,d) := \frac{\ker (d \: A^i \to A^{i+1})}{\mathrm{im}
(d \: A^{i-1} \to A^i)}.$$

If $(B,d_B)$ is another DGA, then a $\mathbb{K}$-linear map $f \: A
\to B$ is called \emph{morphism} if $f(A^i) \subset B^i$, $f$ is
multiplicative, and $d_B \circ f = f \circ d_A$. Obviously, any such
$f$ induces a homomorphism $f^* \: H^*(A,d_A) \to H^*(B,d_B)$. A
morphism of differential graded algebras inducing an isomorphism on
cohomology is called
\emph{quasi-isomorphism}\index{Quasi-Isomorphism}.
\begin{Def}
A DGA $(\mathcal{M},d)$ is said to be \emph{minimal} if
\begin{itemize}
\item[(i)] there is a graded vector space $V = \big( \bigoplus_{i \in
\mathbb{N}_+} V^i \big) = \mathrm{Span} \, \{a_k ~ | ~ k \in I\}$
with homogeneous elements $a_k$, which we call the generators,

\item[(ii)] $\mathcal{M} = \bigwedge V$,

\item[(iii)]  the index set $I$ is well ordered, such that $k < l ~
\Rightarrow |a_k| \le |a_l|$ and the expression for $d(a_k)$ contains
only generators $a_l$ with $l < k$.
\end{itemize}
\end{Def}

We shall say that $(\mathcal{M},d)$ is a \emph{minimal model for a
differential graded algebra}\index{Minimal Model of a! Differential
Graded Algebra} $(A,d_A)$ if $(\mathcal{M},d)$ is minimal and there
is a quasi-isomorphism of differential graded algebras $\rho \: 
(\mathcal{M},d) \to (A,d_A)$, i.e.\ $\rho$ induces an isomorphism 
$\rho^* \: H^*(\mathcal{M},d) \to H^*(A,d_A)$  on cohomology.

The importance of minimal models is reflected by the following
theorem, which is taken from Sullivan's work \cite[Section
5]{Sullivan}.
\begin{Thm}\label{Konstruktion d m M}
 A differential graded algebra $(A,d_A)$ with $H^0(A,d_A) =
\mathbb{K}$ possesses a minimal model. It is unique up to
isomorphism of differential graded algebras.
\end{Thm}

We quote the existence-part of Sullivan's proof, which gives an
explicit construction of the minimal model. Whenever we are going to
construct such a model for a given algebra in this note, we will
do it as we do it in this proof.

\textit{Proof of the existence.} We need the following algebraic
operations to ``add'' resp.\ ``kill'' cohomology.

Let $(\mathcal{M},d)$ be a DGA. We ``add'' cohomology by choosing a
new generator $x$ and setting
$$\widetilde{\mathcal{M}} := \mathcal{M} \otimes \bigwedge(x),~~
\tilde{d}|_{\mathcal{M}} = d, ~~ \tilde{d}(x) = 0,$$ and ``kill'' a
cohomology class $[z] \in H^k(\mathcal{M},d)$ by choosing a new
generator $y$ of degree $k-1$ and setting
$$\widetilde{\mathcal{M}} := \mathcal{M} \otimes \bigwedge(y),~~
\tilde{d}|_{\mathcal{M}} = d, ~~ \tilde{d}(y) = z.$$ Note that $z$
is a polynomial in the generators of $\mathcal{M}$.

Now, let $(A,d_A)$ a DGA with $H^0(A,d_A) = \mathbb{K}$. We set
$\mathcal{M}_0 := \mathbb{K}$,  $d_0 := 0$ and $\rho_0(x) = x$.

Suppose now $\rho_k \: (\mathcal{M}_k,d_k) \to (A,d_A)$ has been
constructed so that $\rho_k$ induces isomorphisms on cohomology in
degrees $\le k$ and a monomorphism in degree $(k+ \nolinebreak 1)$.

``Add'' cohomology in degree $(k+1)$ to get a morphism of
differential graded algebras $\rho_{(k+1),0} \:
(\mathcal{M}_{(k+1),0},d_{(k+1),0}) \to (A,d_A)$ which induces an
isomorphism $\rho_{(k+1),0}^*$ on cohomology in degrees $\le (k+1)$.
Now, we want to make the induced map $\rho_{(k+1),0}^*$ injective on
cohomology in degree $(k+2)$ .

We ``kill'' the kernel on cohomology in degree $(k+2)$ (by
non-closed generators of degree (k+1)) and define $\rho_{(k+1),1} \:
(\mathcal{M}_{(k+1),1},d_{(k+1),1}) \to (A,d_A)$ accordingly. If
there are generators of degree one in
$(\mathcal{M}_{(k+1),0},d_{(k+1),0})$ it is possible that this
killing process generates new kernel on cohomology in degree
$(k+2)$. Therefore, we may have to ``kill'' the kernel in degree
$(k+2)$ repeatedly.

We end up with a morphism $\rho_{(k+1),\infty} \:
(\mathcal{M}_{(k+1),\infty},d_{(k+1),\infty}) \to (A,d_A)$ which
induces isomorphisms on cohomology in degrees $\le (k+1)$ and a
monomorphism in degree $(k+2)$. Now, we are going to set $\rho_{k+1} :=
\rho_{(k+1),\infty}$ and $(\mathcal{M}_{k+1},d_{k+1}) :=
(\mathcal{M}_{(k+1),\infty},d_{(k+1),\infty})$.

Inductively we get the minimal model $\rho \: (\mathcal{M},d) \to
(A,d_A)$. \q
\pagebreak

A \emph{minimal model} $(\mathcal{M}_M,d)$ \emph{of a connected
smooth manifold}\index{Minimal Model of a!Manifold} $M$ is a minimal
model for the de Rahm complex $(\Omega(M),d)$ of differential forms
on $M$. Note that this implies that $(\mathcal{M},d)$ is an algebra
over $\mathbb{R}$. The last theorem implies that every connected
smooth manifold possesses a minimal model which is unique up to
isomorphism of differential graded algebras.

For a certain class of spaces that includes all nilpotent (and hence
all simply-connected) spaces, we can read off the non-torsion part
of the homotopy from the generators of the minimal model. The 
definition of a nilpotent space will be given in the next section.

In general, it is very difficult to calculate the homotopy groups
$\pi_k (X)$ of a given topological space $X$. However, if one is
willing to forget the torsion, with certain assumptions on X, the
rational homotopy groups $\pi_k(X) \otimes \mathbb{Q}$ can be
determined by the theory of minimal models.

In order to relate minimal models to rational homotopy theory, we
need a differential graded algebra over $\mathbb{Q}$ to replace the
de Rahm algebra.

Let $\Delta^n$ be a standard simplex in $\mathbb{R}^{n+1}$ and
$(\Omega_{PL}(\Delta^n),d)$ the restriction to $\Delta^n$ of all
differential forms in $\mathbb{R}^{n+1}$ that can be written as
$\sum P_{i_1 \ldots i_k} dx_{i_1} \ldots dx_{i_k}$, where $P_{i_1
\ldots i_k} \in \mathbb{Q}[x_1 ,\ldots , x_{n+1}]$, together with
multiplication and differential induced by $\mathbb{R}^{n+1}$.

Let $X = \{ (\sigma_i)_{i \in I} \}$ be a path-connected simplicial
complex. Set for $k \in \mathbb{Z}$
$$\Omega^k_{PL}(X) := \{ (\alpha_i)_{i \in I} \,|\, \alpha_i \in
\Omega^k_{PL}(\sigma_i) \wedge ( \sigma_i\subset \partial \sigma_j
\Rightarrow \alpha_j |_{\sigma_i} = \alpha_i) \},$$ and
$\Omega_{PL}(X) := \bigoplus_{k \in \mathbb{Z}} \Omega_{PL}^k(X)$.
It can be verified that the set $\Omega_{PL}(X)$ of so-called
\emph{PL forms}\index{PL Form} is a differential graded algebra over
$\mathbb{Q}$ if we use the multiplication and the differential on
forms componentwise.

\begin{Rem}[{\cite[Remark 1.1.6]{TO}}] 
In fact, PL forms can be defined, along with the minimal model, for 
any CW-complex, say. The process consists of taking the singular 
complex of the space and treating it as a simplicial set amenable to
the PL form construction.
\end{Rem}

For a CW-complex $X$, we define the
\emph{($\mathbb{Q}$-)minimal model}\index{Minimal Model of
a!Simplicial Complex} $\mathcal{M}_{X,\mathbb{Q}}$ of $X$ to be the
minimal model of $(\Omega_{PL}(X),d)$.

\section{Nilpotent spaces}
Already in his paper \cite{Sullivan}, Sullivan shows that for
nilpotent spaces, there is a correspondence between the minimal
model and the rational homotopy. To state this result, we need the
notion of a nilpotent space resp.\ nilpotent module.

Let $G$ be a group, $H$ be a $G$-module, $\Gamma^0_GH := H$ and
$$\Gamma^{i+1}_GH := \langle g.h - h \,|\, g \in G
\wedge h \in \Gamma^i_GH \rangle \subset \Gamma^i_GH$$ for $i \in
\mathbb{N}$.

Then, $H$ is called a \emph{nilpotent
module}\index{Nilpotent!Module} if there is $n_0 \in \mathbb{N}$
such that $\Gamma^{n_0}_GH = \{ 1 \}$.

We recall the natural $\pi_1$-module structure of the higher
homotopy groups $\pi_n$ of a topological space.
For instance, let $(X,x_0)$ be a pointed space with universal cover
$(\widetilde{X},\widetilde{x}_0)$. It is well known that
$\pi_1(X,x_0) \cong D(\widetilde{X})$, the group of deck
transformations of the universal covering. Now, because
$\widetilde{X}$ is simply-connected, every free homotopy class of
self-maps of $\widetilde{X}$ determines uniquely a class of
basepoint preserving self-maps of $\widetilde{X}$ (see e.g.\
\cite[Proposition 4.A.2]{Hatcher}). This means that to every
homotopy class of deck transformations corresponds a homotopy class
of basepoint preserving self-maps (which are, in fact, homotopy
equivalences) $(\widetilde{X},\widetilde{x}_0) \to
(\widetilde{X},\widetilde{x}_0)$. These maps provide induced
automorphisms of homotopy groups
$\pi_n(\widetilde{X},\widetilde{x}_0) \cong \pi_n(X,x_0)$ ($n > 1$)
and this whole process then provides an action of $\pi_1(X,x_0)$ on
$\pi_n(X,x_0)$.

\begin{Def}  \index{Nilpotent!Space}
A path-connected topological space $X$ whose universal covering
exists is called \emph{nilpotent} if for $x_0 \in X$ the fundamental
group $\pi_1(X,x_0)$ is a nilpotent group and the higher homotopy
groups $\pi_n(X,x_0)$ are nilpotent $\pi_1(X,x_0)$-modules for all
$n \in \mathbb{N}$, $n \ge 2$. Note, the definition is independent
of the choice of the base point.
\end{Def}

\begin{Ex} $\,$
\begin{itemize}
\item[(i)] Simply-connected spaces are nilpotent.
\item[(ii)] $S^1$ is nilpotent.
\item[(iii)] The cartesian product of two nilpotent spaces is
nilpotent. Therefore, all tori are nilpotent.
\item[(iv)] The Klein bottle is not nilpotent.
\item[(v)] $P^n(\mathbb{R})$ is nilpotent if and only if $n
\equiv 1(2)$.
\end{itemize}
\end{Ex}

\textit{Proof.} (i) - (iv) are obvious and (v) can be found in
Hilton's book \cite{Hil} on page $165$. \q
\N The main theorem on the rational homotopy of nilpotent spaces is
the following.
\begin{Thm} \label{ratMod}
Let $X$ be a path-connected nilpotent CW-complex with finitely
generated
%fundamental group and finitely generated
homotopy groups. If $\mathcal{M}_{X,\mathbb{Q}} = \bigwedge V$
denotes the $\mathbb{Q}$-minimal model, then for all $k \in \mathbb{N}$ with $k
\ge 2$ holds:
$$\mathrm{Hom}_{\mathbb{Z}}(\pi_k(X),\mathbb{Q}) \cong V^k$$
\end{Thm}

%Erster Satz ist ein Zitat aus Halperin: Lectures on Minimal Models
Using another approach to minimal models (via localisation of spaces
and Postnikow towers), this theorem is proved for example in
\cite{Leh}. The proof that we shall give here is new to the author's
knowledge. We will show the following more general result mentioned
(but not proved) by Halperin in \cite{Halp}.

\begin{Thm}  \label{Satz Halp}
Let $X$ be a path-connected triangulable topological space whose
universal covering exists. Denote by $\mathcal{M}_{X,\mathbb{Q}} =
\bigwedge V$ the $\mathbb{Q}$-minimal model and assume that
\begin{itemize}
\item[(i)] each $\pi_k(X)$ is a finitely generated nilpotent
$\pi_1(X)$-module for $k \ge 2$ and
\item[(ii)] the minimal model for $K(\pi_1(X) , 1)$ has
no generators in degrees greater than one.
\end{itemize}
Then for each $k \ge 2$ there is an isomorphism
$\mathrm{Hom}_{\mathbb{Z}}(\pi_k(X),\mathbb{Q}) \cong V^k$.
\end{Thm}

\begin{Rem}
The homotopy groups of a compact nilpotent smooth manifold are
finitely generated:

By \cite[Satz 7.22]{Hil}, a nilpotent space has finitely generated
homotopy if and only if it has finitely generated homology with
$\mathbb{Z}$-coefficients. The latter is satisfied for compact
spaces. \q
\end{Rem}

The main tool for the proof of the above theorems is a consequence
of the fundamental theorem of Halperin \cite{Halp}. In the next
section, we quote it and use it to prove Theorems \ref{ratMod} and
\ref{Satz Halp}.
\section{The Halperin-Grivel-Thomas theorem}
To state the theorem, let us recall a basic construction for
fibrations.

Let $\pi \: E \to B$ be a fibration with path-connected basis $B$.
Therefore, all fibers $F_b = \pi^{-1}(\{b\})$ are homotopy
equivalent to a fixed fiber $F$ since each path $\gamma$ in $B$
lifts to a homotopy equivalence $L_{\gamma} \: F_{\gamma(0)} \to
F_{\gamma(1)}$ between the fibers over the endpoints of $\gamma$. In
particular, restricting the paths to loops at a basepoint of $B$ we
obtain homotopy equivalences $L_{\gamma} \: F \to F$ for $F$ the
fibre over the basepoint $b_0$. One can show that this induces a
natural $\pi_1(B,b_0)$-module structure on % $H_*(F,\mathbb{Z}$ and on
$H^*(F,\mathbb{Q})$.

\begin{Thm}[{\cite[Theorem 1.4.4]{TO}}]  \label{Fundamentalsatz}
Let $F,E,B$ be path-connected triangulable topological spaces and $F
\to E \to B$ a fibration such that $H^n(F,\mathbb{Q})$ is a
nilpotent $\pi_1(B,b_0)$-module for $n \in \mathbb{N}_+$. The
fibration induces a sequence
$$(\Omega_{PL}(B),d_B) \longrightarrow (\Omega_{PL}(E),d_E)
\longrightarrow (\Omega_{PL}(F),d_F)$$ of differential graded
algebras. Suppose that $H^*(F,\mathbb{Q})$ or $H^*(B,\mathbb{Q})$ is
of finite type.

Then there is a quasi-isomorphism $\Psi \:
(\mathcal{M}_{B,\mathbb{Q}} \otimes \mathcal{M}_{F,\mathbb{Q}} \,,\,
D) \to (\Omega_{PL}(E), d_E)$ making the following diagram
commutative:
%$$\begin{array}{cccccc}
%(\Omega_{PL}(B), d_B) & \longrightarrow & (\Omega_{PL}(E), d_E) &
%\longrightarrow & (\Omega_{PL}(F), d_F) \\ \uparrow \rho_B & &
%\uparrow \Psi & & \uparrow \rho_F \\ (\mathcal{M}_{B,\mathbb{Q}}
%\,,\, D_B) & \longhookrightarrow & (\mathcal{M}_{B,\mathbb{Q}}
%\otimes \mathcal{M}_{F,\mathbb{Q}} \,,\, D) & \longrightarrow &
%(\mathcal{M}_{F,\mathbb{Q}} \,,\, D_F)
%\end{array}$$
$$
\begin{diagram}
(\Omega_{PL}(B), d_B) && \rTo && (\Omega_{PL}(E), d_E) && \rTo &&
(\Omega_{PL}(F), d_F) \\ \uTo^{\rho_B} && && \uTo^{\Psi} && &&
\uTo^{\rho_F} \\ (\mathcal{M}_{B,\mathbb{Q}} \,,\, D_B) && \rInto &&
(\mathcal{M}_{B,\mathbb{Q}} \otimes \mathcal{M}_{F,\mathbb{Q}} \,,\,
D) && \rTo && (\mathcal{M}_{F,\mathbb{Q}} \,,\, D_F) \\
\end{diagram}
$$
Furthermore,  the left and the right vertical arrows are the minimal
models. Moreover, if $\mathcal{M}_F = \bigwedge V_F$, there is an
ordered basis $\{v_i^F \,|\, i \in I \}$ of $V_F$ such that for all
$i,j \in I$ holds $D(v_i^F) \in \mathcal{M}_B \otimes
(\mathcal{M}_F)_{<v_i^F}$ and $(v_i^F < v_j^F \Rightarrow |v_i^F|
\le |v_j^F|)$. \q
\end{Thm}

\begin{Rem}
In general, $(\mathcal{M}_{B,\mathbb{Q}} \otimes
\mathcal{M}_{F,\mathbb{Q}} \,,\, D)$ is not a minimal differential
graded algebra and $D|_{\mathcal{M}_{F,\mathbb{Q}}}\ne D_F$ is
possible.
\end{Rem}

We need some further preparations for the proofs of the above
theorems. The first is a reformulation of the results $3.8 - 3.10$
in \cite{Hil}. It justifies the statement of the next theorem.

\begin{Prop}\label{Torsion nil}
Let $G$ be a finitely generated nilpotent group. Then the set $T(G)$
of torsion elements of $G$ is a finite normal subgroup of $G$ and $G
/ T(G)$ is finitely generated. \q
\end{Prop}

\begin{Thm}  \label{Knil}
Let $G$ be a finite generated nilpotent group and denote by $T(G)$
its finite normal torsion group.

Then $K(G,1)$ and $K(G/T(G) , 1)$ share their minimal model.
\end{Thm}

\textit{Proof.} Since $T(G)$ is finite and $\mathbb{Q}$ is a field,
we get from \cite[Section 4.2]{Evens} $H^n(K(T(G) , 1)), \mathbb{Q})
= \{0\}$ for $n \in \mathbb{N}_+$. The construction of the minimal
model in the proof of Theorem \ref{Konstruktion d m M} implies that
$\mathcal{M}_{K(T(G) , 1),\mathbb{Q}}$ has no generators of degree
greater than zero. Now, the theorem follows from the preceding one,
applied to the fibration $K(T(G) , 1) \to K(G,1) \to K(G/T(G) , 1)$.
\q

\begin{Lemma}
Let $X$ be topological space with universal covering $\mathfrak{p}
\: \widetilde{X} \to X$.

Then, up to weak homotopy equivalence of the total space, there is a
fibration $\widetilde{X} \to X \to K := K(\pi_1(X),1)$. Moreover,
for a class $[\gamma] \in \pi_1(K) \cong \pi_1(X)$ the homotopy
equivalences $L_{[\gamma]} \: \widetilde{X} \to \widetilde{X}$
described at the beginning of this section are given by the
corresponding deck transformations of $\mathfrak{p}$.
\end{Lemma}

\textit{Proof.} Denote by $\pi \: E \to K(\pi_1(X),1)$ the universal
principal $\pi_1(X)$-bundle. Regard on $E \times \widetilde{X}$ the
diagonal $\pi_1(X)$-action. Then, the fibre bundle $$\widetilde{X}
\longrightarrow \big( (E \times \widetilde{X}) / \pi_1(X) \big)
\longrightarrow K$$ has the desired properties. \q

\subsection{Proof of Theorem \ref{Satz Halp}:} Let $X$ be as in
the statement of the theorem. For simply-connected spaces, the
theorem was proven in \cite[Theorem 15.11]{FHT}. Now, the idea is to
use this result and to consider the universal cover $\mathfrak{p} \:
\widetilde{X} \to X$. Denote by $\mathcal{M}_{\widetilde{X},
\mathbb{Q}} = \bigwedge \widetilde{V}$ and $\mathcal{M}_{X,
\mathbb{Q}} = \bigwedge V$ the $\mathbb{Q}$-minimal models. We shall 
show
\begin{equation} \label{1}
\forall_{k \ge 2}~ V^k \cong \widetilde{V}^k.
\end{equation}
This and the truth of the theorem for simply-connected spaces
implies then the general case $$\forall_{k \ge 2}~ V^k \cong
\widetilde{V}^k \cong \mathrm{Hom}_{\mathbb{Z}}(\pi_k(\widetilde{X})
, \mathbb{Q}) = \mathrm{Hom}_{\mathbb{Z}}(\pi_k(X) , \mathbb{Q}).$$

It remains to show (\ref{1}): Since $X$ is triangulable, $X$ and
$\widetilde{X}$ can be seen as CW-complexes. Therefore, up to weak
homotopy, there is the following fibration of CW-complexes
$$\widetilde{X} \longrightarrow X \stackrel{\pi}{\longrightarrow}
K(\pi_1(X),1) =: K.$$

We prove below:
\begin{eqnarray}
\label{2} & H^*(\widetilde{X},\mathbb{Q}) \mbox{ is of finite type.} & \\
\label{3} & H^*(\widetilde{X},\mathbb{Q}) \mbox{ is a nilpotent }
\pi_1(X)\mbox{-module.} &
\end{eqnarray}
Then Theorem \ref{Fundamentalsatz} implies the existence of a
quasi-isomorphism $\rho
%\: (\mathcal{M}_{K,\mathbb{Q}} \otimes \mathcal{M}_{\widetilde{X}, \mathbb{Q}}
%\,,\, D) \to (\Omega_{PL}(X),d_X)
$ such that the following diagram commutes:
%$$\begin{array}{cccccc}
%(\Omega_{PL}(K), d_K) & \longrightarrow & (\Omega_{PL}(X), d_X) &
%\longrightarrow & (\Omega_{PL}(\widetilde{X}), d_{\widetilde{X}}) \\
%\uparrow \rho_K & & \uparrow \rho & & \uparrow \rho_{\widetilde{X}} \\
%(\mathcal{M}_{K,\mathbb{Q}} \,,\, D_K) & \longhookrightarrow &
%(\mathcal{M}_{K,\mathbb{Q}} \otimes
%\mathcal{M}_{\widetilde{X},\mathbb{Q}} \,,\, D) & \longrightarrow &
%(\mathcal{M}_{\widetilde{X},\mathbb{Q}} \,,\, D_{\widetilde{X}})
%\end{array}$$
$$
\begin{diagram}
(\Omega_{PL}(K), d_K) &&  \rTo && (\Omega_{PL}(X), d_X) &&
\rTo && (\Omega_{PL}(\widetilde{X}), d_{\widetilde{X}}) \\
\uTo^{\rho_K} && && \uTo^{\rho} && && \uTo^{\rho_{\widetilde{X}}} \\
(\mathcal{M}_{K,\mathbb{Q}} \,,\, D_K) && \rInto &&
(\mathcal{M}_{K,\mathbb{Q}} \otimes
\mathcal{M}_{\widetilde{X},\mathbb{Q}} \,,\, D) && \rTo &&
(\mathcal{M}_{\widetilde{X},\mathbb{Q}} \,,\, D_{\widetilde{X}})
\end{diagram}
$$
Finally, we shall see
\begin{equation} \label{4}
(\mathcal{M}_{K,\mathbb{Q}} \otimes
\mathcal{M}_{\widetilde{X},\mathbb{Q}} \,,\, D) \mbox{ is a minimal
differential graded algebra}
\end{equation}
and this implies (\ref{1}) since $\mathcal{M}_K$ has no generators
of degree greater than one by assumption (ii).

We still have to prove (\ref{2}) - (\ref{4}):

By assumption (i), $\pi_k(X) = \pi_k(\widetilde{X})$ is finitely
generated for $k \ge 2$. Since simply-connected spaces are
nilpotent, \cite[Satz 7.22]{Hil} implies the finite generation of
$H_*(\widetilde{X},\mathbb{Z})$ and (\ref{2}) follows.

(\ref{3}) is the statement of Theorem 2.1 $(i) \Rightarrow (ii)$ in
\cite{Hil76} -- applied to the action of $\pi_1(X)$ on
$\pi_i(\widetilde{X})$.

ad (\ref{4}): By assumption (ii), $\mathcal{M}_K$ has no generators
in degrees greater than one, i.e.\ $\mathcal{M}_{K,\mathbb{Q}} =
\bigwedge \{v_i \,|\, i \in I \}$ with $|v_i| = 1$. The construction
of the minimal model in the proof of Theorem \ref{Konstruktion d m
M} implies that the minimal model of a simply-connected space has no
generators in degree one, i.e.\
$\mathcal{M}_{\widetilde{X},\mathbb{Q}} = \bigwedge \{w_j \,|\, j
\in J \}$ with $|w_j| > 1$. We expand the well orderings of $I$ and
$J$ to a well ordering of their union by $\forall_{i \in I} ~ \forall_{j \in J}~~ i<j$.
 Theorem \ref{Fundamentalsatz} implies that
$D(w_j)$ contains only generators which are ordered before $w_j$.
Trivially, $D(v_i)$ also has this property, so we have shown
(\ref{4}) and the theorem is proved. \q

\subsection{Proof of Theorem \ref{ratMod}:}
Let $X$ be a path-connected nilpotent CW-complex with finitely
generated fundamental group and finitely generated homotopy. By
Theorem \ref{Satz Halp}, we have to show that the minimal model of
$K(\pi_1(X),1)$ has no generators in degrees greater than one.

Theorem \ref{Knil} implies that it suffices to show that $K(\pi_1(X)
/ T ,1 )$ has this property, where $T$ denotes the torsion group of
$\pi_1(X)$. $\Gamma := \pi_1(X) / T$ is a finitely generated
nilpotent group without torsion. By \cite[Theorem 2.18]{Rag},
$\Gamma$ can be embedded as a lattice in a connected and
simply-connected nilpotent Lie group $G$. Therefore, the nilmanifold
$G / \Gamma$ is a $K(\Gamma,1)$ and from \cite[Theorem 3.11]{ipse}
follows that its minimal model has no generators in degrees greater
than one. \q

\begin{Ac}
The results presented in this paper are parts of my dissertation that I wrote 
under the supervision of H.\ Geiges. 
I wish to express my sincerest gratitude for his support.
Moreover, I wish to thank St.\ Halperin. I have profited from his suggestions.
\end{Ac}

\begin{Rem}
St.\ Halperin told me that smooth triangulations are unnecessary to do rational
homotopy for smooth differential forms, as is already presented in his works 
\cite{Halp} and \cite{FHT}.
\end{Rem}

\noindent
\textsc{Christoph Bock\\ Department Mathematik\\ Universität Erlangen-N\"urn\-berg\\ Cauerstraße 11\\ 91058 Erlangen\\ Germany}

\noindent
\textit{e-mail:} \verb"bock@mi.uni-erlangen.de"
\end{document}